# Loss of Tension in an Infinite Membrane with Holes Distributed by Poisson Law

M. V. Menshikov[*]

menchik@math.usp.ac.br

Department of Mathematics and Statistics, University of San Paulo
and Faculty of Mechanics and Mathematics, Moscow State University

K. A. Rybnikov[†][‡]

rybnikov@mast.queensu.ca

Department of Mathematics and Statistics, Queen's University and
Fields Institute for Research in Mathematical Sciences

S. E. Volkov,

svolkov@fields.utoronto.ca

Fields Institute for Research in Mathematical Sciences
and Department of Mathematics and Statistics, York University

August 12, 1999

## Abstract

If one randomly punches holes in an infinite tensed membrane, when does the tension cease to exist? This problem was introduced by R. Connelly in connection with applications of rigidity theory to natural sciences. We outline a mathematical theory of tension based on graph rigidity theory and introduce several probabilistic models for this problem. We show that if the "centers" of the holes are distributed in $\mathbb{R}^2$ according to Poisson law with parameter $\lambda > 0$, and the distribution of sizes of the holes is independent of the distribution of their centers, the tension vanishes on all of $\mathbb{R}^2$ for any value of $\lambda$. In fact, it follows from a more general result on the behavior of iterative convex hulls of connected subsets of $\mathbb{R}^d$, when the initial configuration of subsets is distributed according to Poisson law and the sizes of the elements of the original configuration are independent of this Poisson distribution. For the latter problem we establish the existence of a critical threshold in terms of the number of iterative convex hull operations required for covering all of $\mathbb{R}^d$. The processes described in the paper are somewhat related to bootstrap and rigidity percolation models.

---

[*]most of the results were developed while the author was visiting the Fields Institute

[†]the work of K. A. Rybnikov was supported in part by fellowships from the Fields Institute

[‡]corresponding author

# 1 Introduction

Let $M$ be a tensed membrane (film) clamped on the boundary. If one makes a small convex hole in this membrane, the tension redistributes over the rest of the membrane. However the situation is more complicated when there are several hole in the membrane. It is rather clear that if you have a non-convex hole or a couple of convex overlapping holes, tension ought to vanish on the convex hull of this set (see **Fig. 1.1**).

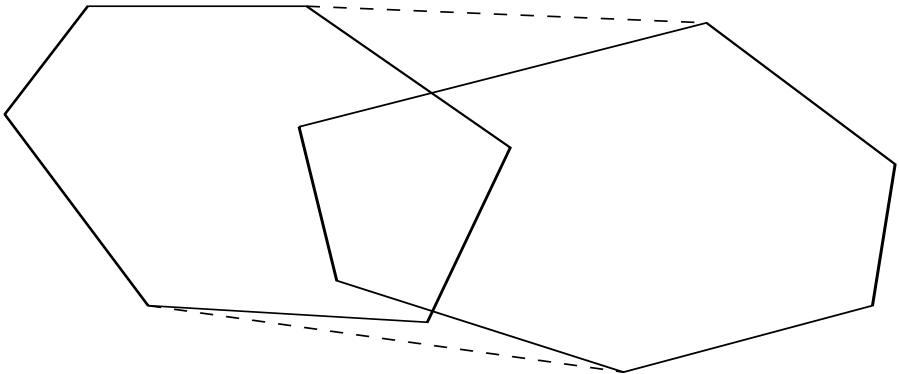

**Figure 1.1**

However, it is less intuitive that tension may vanish at some subset of the complement of a collection of convex non-overlapping holes.

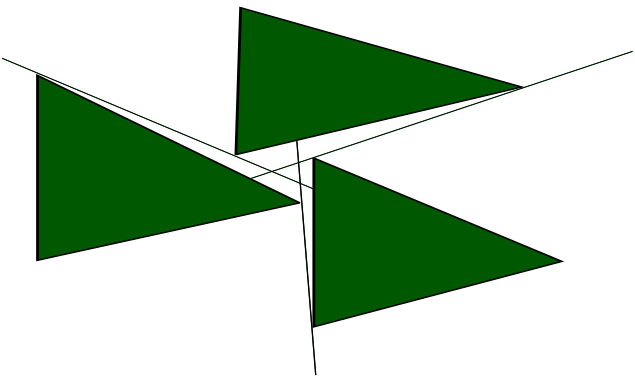

**Figure 1.2**

For example, the convex hull of three holes shown on **Fig. 1.2** cannot support tension which may be verified in practice with a sheet of some elastic material and scissors. Therefore, if the area where tension ought to vanish is interpreted as defective, all three polygons on **Fig. 1.1** ought to coalesce into one big defect. Such *coalescence effect of a "pinwheel configuration"* was first noticed by Connelly (1998). A mathematical explanation of this fact will be given in Subsection 1.2.

In the following subsection we give a formal treatment of tension in the language of graph rigidity. Our mathematical treatment of a tensed medium is discrete in nature and based on the ideas from rigidity theory, as opposed to classical methods of mechanics of continuous media that use smooth functions, tensor fields, and differential equations for describing elastic properties of materials. Notice that in the planar case holes with smooth boundary can be well approximated by polygons and that the theory based on graph rigidity seem is better suited for computer simulation. The existence of a "discrete" tension in the complement of the holes has some important implications for rigidity properties of large planar graphs realized in the complement of the holes. This graph theoretic treatment of tension was introduced in (Connelly, Mitchel, Rybnikov (1998, 1999)). In this paper we are primarily interested in the behavior of the system with very large or infinite number of holes. More specifically, we look at the question of existence of a subset supporting tension in the complement of an infinite set of holes in $\mathbb{R}^2$. It follows directly from definitions (see Subsection 1.1) that if such set exists, it is infinite and a.s. unique, which is, apparently in good agreement with intuition. We consider a few models for the distribution of holes in an infinite membrane (Section 2), and completely resolve the question of the existence of a tensed subset for a "naïve" Poisson model suggested by Connelly. In this model holes are distributed according to a Poisson law, and their shapes are IID. In fact we prove a somewhat stronger statement, which we hope may have some reasonable interpretation (see Subsection 1.1). The percolation models that we describe in this paper are closely related to so-called "bootstrap percolation" introduced on trees by Chalupa, Leath and Reich (1979) and later on $d$-dimensional lattices by Kogut and Leath (1981). In these models, points are independently occupied with a low density and the resulting configuration is taken as the initial state for dynamics based on some collection of local rules, in which the occupation status of a point is updated according to the configuration of its neighbors. The rigorous analysis of this model have been conducted by Aizenman and Lebowitz (1988). Another work, showing that the critical probability $p_c = 1$ for some bootstrap percolation models was by van Enter (1987). Since then, there has been many publications devoted to this topic. For the latest results on bootstrap percolation see Dehghanpour and Schonmann (1997) and Gravner and McDonald (1997). Gravner and Schonmann (1999, personal communication) had some ideas how to extend their results for continuous case, but as far as the authors know, there has been no published results on the model which we describe and analyze below.

## 1.1 Rigidity and its Applications

Over the past two decades 2-dimensional and 3-dimensional random central-force networks have been used by physicists for modeling the elastic behavior of glasses within the framework of effective medium theory (see Duxbury, Thorpe, Jacobs, Moukarzel (1983, 1995, 1996, 1998, 1999)). It turns out that real glasses are well represented by generic random networks. The success of these methods resulted in good characterization of elastic properties of glasses like $Ge_xAs_ySe_{1-y}$ (Thorpe (1983)). The rigidity analysis of central-force networks has been also used for characterization of physical properties of other substances such has proteins and semiconductors (see Thorpe and Duxbury (1999)).

While discussing the property of a membrane to be able to support tension, we keep in mind

(large) networks in the complement of the holes. It turns out that the existence of tension in the compliment of the holes implies some important properties for a network realized in the complement. At this point we need to introduce some mathematical terminology. Notice that in mathematics of rigidity there is a tendency to use term *framework* instead of *network* preferred by physicists.

A bar-and-joint framework is a realization of a graph in $\mathbb{R}^d$. Denote by $F(E;V,V_0)$ a framework in $\mathbb{R}^d$ with the edge set $E$, and the vertex set $V$, with pinned (fixed in $\mathbb{R}^d$) subset of vertices $V_0 \subset V$. Denote by $\mathbf{v}_i$ the vector of coordinates of vertex $v_i \in V$.

**Definition 1.1** *An equilibrium stress (or self-stress) is an assignment of real numbers $s_{ij} = s_{ji}$ to the edges, a tension if the sign is positive, or a compression if the sign is negative, so that the equilibrium conditions*

$$\sum_{(ij) \in E} s_{ij}(\mathbf{v}_j - \mathbf{v}_i) = \mathbf{0}$$

*hold at each vertex $\mathbf{v}_i \in V \backslash V_0$ (see* **Fig. 1.3***).*

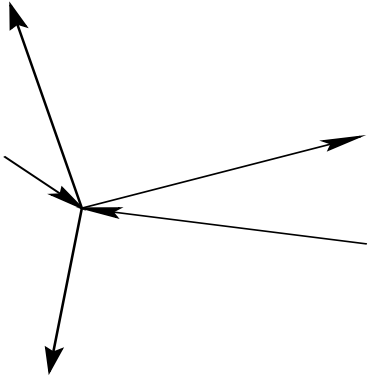

**Figure 1.3**

A finite framework $F(E;V)$ is called *rigid* in $\mathbb{R}^d$ if and only if there is a neighborhood $N(F)$ of $F$ in the space of parameters $\mathbb{R}^{dV}$ such that any other $\mathbb{R}^d$ realization of graph $(E,V)$ from $N(F)$ with the same lengths of all edges is congruent to $F(E;V)$. Sometimes it is interesting to study rigidity of graphs some of whose vertices are pinned. Framework $F(E;V,V_0)$ is called rigid with pinned vertices $V_0$ if $F$ has a neighborhood $N(F)$ in $\mathbb{R}^{dV}$ such that any other realization of $(E;V)$ from $N(F)$ with the same lengths of all edges and the set of pinned vertices inherited from the set of pinned vertices of $F$ is congruent to $F$. If the "neighborhood condition" is dropped the framework is called globally rigid. A framework that is not rigid is called *flexible*. A framework $F(E;V,V_0)$ that supports an all non-zero equilibrium tension is referred to as a *spider web*.

**Definition 1.2** *An framework $F(E;V,V_0)$ (possibly infinite) is referred to as rigid if any finite subgraph of $F(E;V,V_0)$ is contained in some rigid finite subgraph of $F(E;V,V_0)$.*

By polygonal partition of a planar set with piecewise-linear or no boundary we mean an edge-to-edge partition of this set into convex bounded polygons. The 1−skeleton of a partition is a framework whose vertex set is the vertex set of the partition, and whose edge set is the partition edge set. A polygonal partition which consists of triangles is called a *triangulation*.

**Definition 1.3** *Let $M$ be a set with polygonal boundary in $\mathbb{R}^2$ ($M$ might be all of $\mathbb{R}^2$), and let $\mathcal{H}$ be a collection of open polygons in $M$, such that the number of polygons intersecting any compact subset of $\mathbb{R}^2$ is finite. We call the elements of $\mathcal{H}$* holes *and denote by $H$ the pointwise union of the holes. We say that $M \backslash H$* supports tension *if $M \backslash H$ admits a partition with the set of edges $E$ and vertices $V$, such that the framework $(E, V, V \cap \partial M)$ is a spider web.*

Evidently, a polygonal partition can be replaced by a triangulation in this definition without loss of generality. A direct generalization of this definition to the case of general dimension is possible, but not quite natural, since not all spider webs in dimensions higher than 2 can be interpreted as 1-skeletons of polyhedral partitions (see Connelly and Whiteley (1993, 1996)). In the planar case the situation is simplified by the fact that any spider web with self-intersections can be turned into the 1-skeleton of a polygonal partition by adding points of self-intersections to the vertex set of the framework, and modifying the edge set accordingly. The cone of tension of the new partition contains the cone of tension of the original one. A somewhat more natural definition would the one in which the existence of some spider web in the complement of the holes is required. Let us now make some observations about holes. First, if a hole is non-convex, then there is no triangulation of the complement such that its 1-skeleton minus the boundary of $M$ supports equilibrium tension. For instance, the equilibrium of forces at vertex **v** on **Fig. 1.4** is impossible if all edges incident to this vertex are under tension. Therefore, if two holes overlap, and their union is not convex, tension vanishes on all of their convex hull. This is called the *coalescence effect of overlapping holes*.

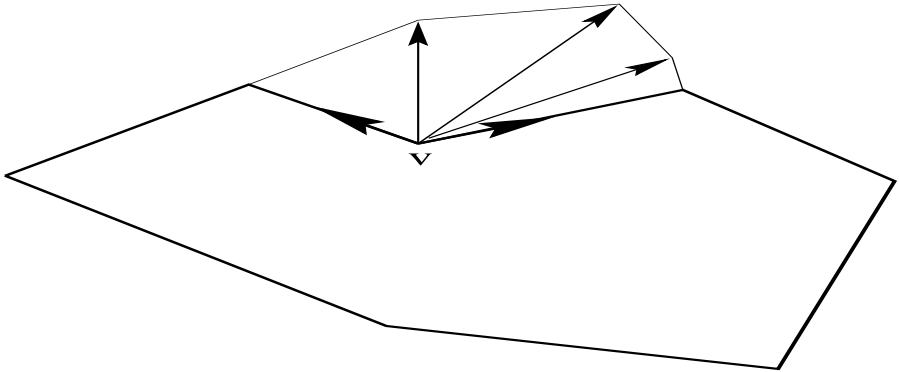

**Figure 1.4**

Notice that the above observation about the convexity of connected components of the complement of the set supporting tension is not valid for dimensions greater than three, whatever definition of tension we adopt. For instance, a three-dimensional polytonal hole can have a saddle point, to which a number of edges lying in the complement of the holes can be attached so that this point will be in static equilibrium. However, a set supporting tension in $\mathbb{R}^d$ cannot have points of strict convexity.

It is natural to call subsets of $M$ that cannot support tension *defects*. The primary defects are, obviously, holes themselves. The following definitions iteratively define the notion of a defect of $k^{th}$ generation. Note that in these definitions the polygonality of the holes is not important.

**Definition 1.4** *Let $\mathcal{H}$ be a set of holes. Elements of $\mathcal{H}$ are called defects of $0^{th}$ generation.*

**Definition 1.5** *A connectivity component (understood topologically) of defects of $k^{th}$ generation is referred to as a $k-$cluster.*

**Definition 1.6** *A defect of $(k+1)^{th}$ generation is the convex hull of a $k-$cluster.*

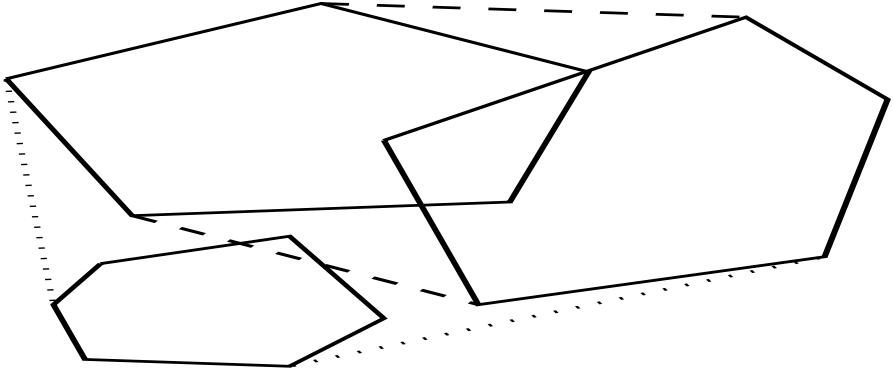

**Figure 1.5**

**Fig. 1.5** provides an illustration to our definitions. Here we used the solid line for the boundary of holes ($\equiv$ defects of $0^{th}$ generation), the dash line for the boundary of a defect of the first generation, and dotted line for the boundary of a defect of the second generation. Note that our definition of defect does not cover all subsets that cannot support tension, since an area that fails to support tension need not belong to the convex hulls of a connected component of the set of holes as it was noticed above (see **Fig. 1.2** and explanations in Section 1.2). However, the effect of a "pinwheel configuration" is insignificant for the "naïve" Poisson percolation model that we analyze in this paper. Even without taking this effect into account, we show that tension disappears on all of $\mathbb{R}^2$.

The existence of a triangulation which supports tension may have some interesting implications for the modeling of physical properties of materials with networks of Hooke springs and geometry of convex surfaces (see Subsection 1.2). Among other things it implies that every triangulation of $M\backslash H$ is *globally* rigid not only in $\mathbb{R}^2$, but also in $\mathbb{R}^3$, provided all its vertices lying on the boundary of $M$ are pinned (Connelly (1993)).

The elasticity (rigidity) properties of a glass are related to how amenable the glass is to continuous deformations which require little energy. From a physical point of view, it is not enough to declare that the distance constraints force the structure to have only one configuration, since the bonds in a physical network do not behave as ideal edges in a framework. There should be a way of describing the behavior of the system as it is perturbed. That is why physicists often consider the energy function defined on the edges of a network of Hooke springs, i.e. each spring has some optimal length at which its energy is minimal, stretching or shortening a string increases the energy of this connection

A tensegrity framework is a generalization of this model where besides Hooke springs there are members whose energy increases with the distance, and members whose energy decreases with the distance. In context of energy considerations it is often useful to work with the notion of tensegrity framework (Connelly, Whiteley (1996)). In a tensegrity framework all edges are partitioned into three types, cables $E_+$, struts $E_-$, and bars $E_0$, i.e. $E = E_0 \cup E_+ \cup E_-$. Together, struts, cables and bars are called members. Let $F(E, E_0, E_+, E_-; V, V_0)$ be some tensegrity framework in $\mathbb{R}^d$, and denote by $l_{ij}^0$ the length of a member $(ij)$. If a cable is stretched the energy in the cable increases. If a strut is shortened, the energy in it decreases. Any change in the length of a bar gets the energy to increase. Therefore networks of Hooke springs are bar tensegrities from mathematical point of view. Thus the energy $\mathfrak{H}_{ij}$ of member $(ij)$ considered as the function of its squared length $l_{ij}$ is monotone increasing if $(ij)$ is a cable,
is monotone decreasing if $(ij)$ is a strut,
has a strict local minimum at $l_{ij}^0$ called the equilibrium length.

It is natural to define the energy function of a finite tensegrity framework as the sum of the energy functions of its members. Thus

$$\mathfrak{H} = \frac{1}{2} \sum_{(ij) \in E} \mathfrak{H}_{ij}(|\mathbf{v}_j - \mathbf{v}_i|^2). \quad (1)$$

The simplest way to define the energy function when all members are bars is as follows

$$\mathfrak{H} = \frac{1}{2} \sum_{(ij)} a_{ij} n_{ij} (l_{ij} - l_{ij}^0)^2, \quad (2)$$

where the sum is over all ordered pairs of vertices of the framework, $n_{ij}$ is 0 or 1 depending on whether there is a bond between $i$ and $j$, $l_{ij}$ is the length of the bond between $i$ and $j$, $l_{ij}^0$ is the equilibrium bond length, and $a_{ij} > 0$ is the spring constant of the bond $(ij)$.

In the spirit of the definition of equilibrium stress we assume that a strut can support only compression, a cable can support only tension, and a bar can be under both types of stress. For more detailed information on tensegrities see the works of Connelly and Whiteley (1995), and Connelly (1993).

Consider now the energy $\mathfrak{H}$ as a function of the *coordinates of the vertices* of the framework. A finite tensegrity framework $F$ in $\mathbb{R}^d$ with pinned vertices $V_0 \subset V$ is called *prestress stable* if

1) the differential of $\mathfrak{H}$ (considered as a function on the space of parameters $\mathbb{R}^d V$) at the point corresponding to $F$ is zero,

2) the second differential of $\mathfrak{H}$ is a positive semidefinite quadratic form whose kernel restricted to infinitesimal motions leaving $V_0$ unmoved consists of trivial infinitesimal motions of the framework.

**Definition 1.7** *An infinite tensegrity framework is called prestress stable if every its finite subgraph $G$ is contained in a prestress stable subgraph whose un-pinned vertices contain the vertex set of $G$.*

The concept of prestress stability is due to Connelly and Whiteley (1993, 1996) and basically accounts for local minima of the energy function. It turns out that if the complement of $H$ in $M$ supports tensions, then for any triangulation $T$ of $M \setminus H$ there is a an interpretation of the edges of $T$ as struts and cable such that the resulting tensegrity is prestress stable in 3-space with $\partial M$ pinned. In particular it implies that $T$ is prestress stable as a bar framework in 3-space with $\partial M$ pinned.

So far we discussed only statics of framework. There is an interesting implication of the non-existence of tension in the complement of a collection of holes for mechanics. It was conjectured by Connelly (1998) that in this case the complement has a triangulation which is a flexible in 3-space with boundary pinned. There also are interesting connections between our problem and convex geometry that in its original form are due to Maxwell (1864, 1869-1872) and Cremona (1872). They are outlined in Section 1.2.

Let us summarize the implications of the existence of tension in the complement of the holes (see also Connelly, Mitchel, and Rybnikov (1998, 1999)).

**Proposition 1.8** *Let $M$ be a convex set of $\mathbb{R}^2$ with polygonal or no boundary, and let $\mathcal{H}$ be a collection of convex open polygons in $M$, possibly overlapping. Suppose $M \backslash (H \cup \partial M)$ supports tension. Let $P$ be a convex polygonal partition of $M \backslash H$ with the set of vertices $V$ and the set of edges $E$, and let $F(E; V, V \cap \partial M)$ be the corresponding framework with pinned vertices $V \cap \partial M$. Then*

*1) framework $F(E; V, V \cap \partial M)$ is globally rigid in $\mathbb{R}^3$ as a bar framework;*

*2) there is a partition of $E$ into two sets $E_+$ and $E_-$, such that if members of $E_+$ are cable and members of $E_-$ are struts, the resulting tensegrity is prestress stable in $\mathbb{R}^3$.*

**Conjecture 1.9** *(Connelly) Let $M$ be a convex subset of $\mathbb{R}^2$ with polygonal or no boundary, and let $\mathcal{H}$ be a collection of convex open polygons in $M$. Suppose $M \backslash (H \cup \partial M)$ does not support tension. Then there is a triangulation of $M \backslash H$ which is flexible in 3-space with the boundary of $M$ pinned.*

We say that tension vanishes at point $p \in M$ if there is no subset of $M \backslash H$ supporting tension that covers this point. In the following subsection we show that for a finite system $\mathcal{H}$ any subset of $M \backslash H$ which supports tension is, in fact, contained in the maximum subset $S_{max} \subseteq M \backslash H$ supporting tension. Thus, when the number of holes is finite, $M \backslash H$ can be partitioned into two polygonal subsets, the maximum subset supporting tension and its complement where tension vanishes. We conjecture that this is also true for infinite case ( see Section 1.2). The above proposition suggests an interpretation of the area that supports tension as stable one, and of the area that cannot support tension as amenable to deformations requiring little energy.

One of the applications we have in mind is percolation of defects in presence of crystallization. Here $\mathcal{H}$ is not a set of holes but rather a dissemination of some alien substance that have poor contact with our material. Crystallization cannot take place in a region where the tension has been lost, but we assume that the system has some memory. The reason why we assume that crystallization cannot occur when tension has been lost is that the loss of tension means the loss of a rigid structure as it is suggested by Proposition 1.8 and Conjecture 1.9. We do not claim that the rigidity of flexibility of a triangulation of the complement of the defects directly explains the behavior of an actual physical system. However, if Connelly's Conjecture 1.9 is correct, the existence of tension is a *necessary and sufficient condition of the global rigidity of all triangulations of a planar set with polygonal boundary.*

We assume that tension cannot be lost instantaneously, since there is some contact between the material and the alien substance. Before tension disappears on defects of $k^{th}$ generation it disappears on defects of $(k-1)^{th}$ generation. In section 2 we show that there is $N \in \mathbb{Z}_+$ such that with probability 1 every point of $\mathbb{R}^2$ is covered by a defect of $N^{th}$ generation. It can be interpreted as that there is time $t_c(\mathcal{H})$ such that if the system has not crystallized by time $t_c(\mathcal{H})$ it will never crystallize.

## 1.2 Tension and convexity

J. C. Maxwell (1864, 1869-1872) and later L. Cremona (1872) observed that a framework in $\mathbb{R}^2$ which can be interpreted as the 1-skeleton of a 2-sphere admits an all non-zero equilibrium stress if and only if this framework is the projection of a polyhedral sphere in $\mathbb{R}^3$. This is also true for a partition of $\mathbb{R}^2$ into convex polygons, but in this case an infinite polyhedral surface plays the role of a 2-sphere. These theorems were later rigorously proved by Crapo and Whiteley (1982, 1984, 1993). General theory of the relationships between stresses, liftings, and Voronoi diagrams can be found in Rybnikov (1999). In fact, this connection between stresses and convex surfaces provides much insight into the problem of understanding the redistribution of stress caused by punched holes. It turns out that in Maxwell construction a *tension* on the 1-skeleton of a polygonal partition corresponds to a *convex lifting* of the partition. Therefore all of the complement of the holes supports tension if and only if the holes can be lifted to *different* faces of some convex surface. Projecting such surface back on the plane we get a partition of the complement of the holes whose skeleton support an all non-zero tension.

**Proposition 1.10** *Let $\mathcal{H}$ be a discrete system of holes in $\mathbb{R}^2$. The complement of $H$ supports a tension if and only if there is a continuous function defined $\mathbb{R}^2$ whose graph is a convex polyhedral surface in $\mathbb{R}^3$ such that all connected components of $H$ are vertical projections of different facets of this surface*

For example, on

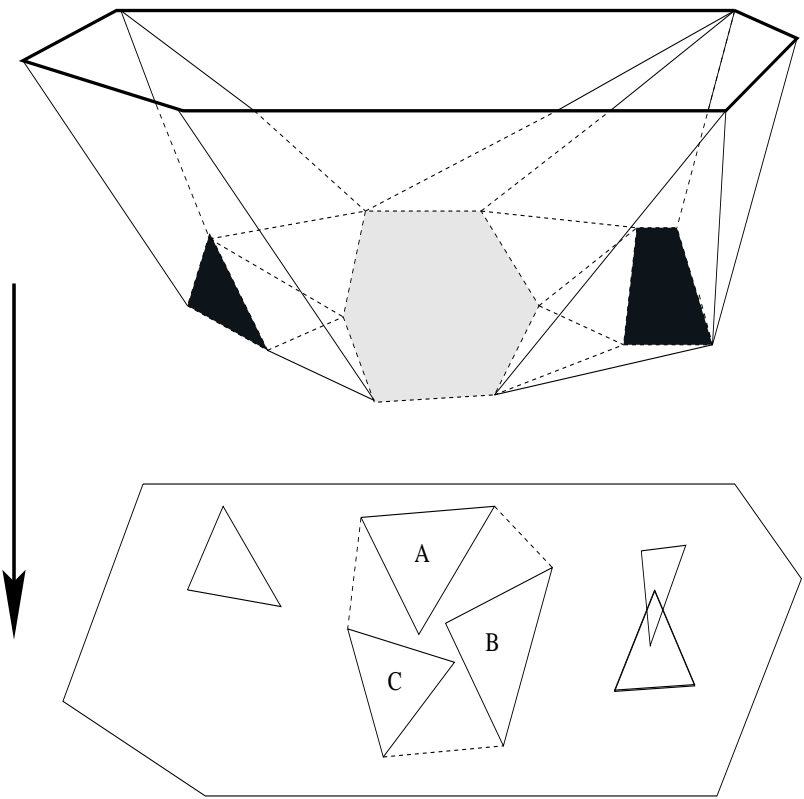

**Figure 1.6**

For example, on **Fig. 1.2** we see that no convex piecewise-linear function can lift the triangles to three different facets, since it would imply that these triangles can be mutually separated in $\mathbb{R}^2$ by three lines meeting at one point. This is clearly impossible for the configuration depicted on **Fig. 1.2**. Similarly, on **Fig. 1.6** we see that the complement of the original system of holes (thin solid line) does not support a tension. However, the system of enlarged holes (dash line) supports a tension. Holes $A$, $B$, and $C$ form a "pinwheel configuration" like one on **Fig. 1.2**.

As we have seen the Maxwell correspondence helps to understand the redistribution of tension in the case of a finite number of holes. In fact it also suggest an algorithm for computing the maximum subset of the complement that supports a tension (Connelly, Mitchel, Rybnikov (1998, 1999)). The following theorem shows that one can exactly define the maximum subset of the complement that supports a tension, and illustrates the use of Maxwell correspondence.

**Lemma 1.11** *Let $M$ be a convex polygon in $\mathbb{R}^2$ ($M$ might be all of $\mathbb{R}^2$), and let $\mathcal{H}$ be a finite collection of convex polygons in $M$. There is a set $H_{min}$ such that*

*1) $H \subset H_{min}$,*

*2) $M \backslash H_{min}$ supports tension,*

*3) There is no proper subset of $H_{min}$ satisfying the first two conditions.*

*Furthermore, $H_{min}$ is unique.*

*Proof.* Let $H'$ be a set of holes that satisfies the first two conditions, namely, $H \subset H'$ and $M \backslash H'$ supports tension. Recall, that any system of holes satisfying the first two conditions is actually a collection of edge-disjoint convex polygons, i.e. no pair of polygons in the collection share an edge.

Suppose not all of the vertices of $H'$ come from the vertices of $H$. We show that in this case $H'$ is not minimal, which means that there is a subset of $H'$ that satisfies Conditions 1 and 2 and is contained in $H'$. Let $D$ be a polygonal hole of $H'$ with a vertex $v$ which does not belong to the vertex set of $H$. Vertex $v$ can be separated from all other vertices of $D$ by a line passing though some two vertices of $H$. The resulting triangle $\Delta$ can be cut off $D$ without violating the property of $M \backslash (H' \backslash \Delta)$ to support tension. Let us illustrate it with the use of Maxwell correspondence. $M \backslash H'$ has a triangulation that lifts to a convex surface. Denote the lifting map by $L$. Let us modify the surface by cutting off the corner of the vertex $L(v)$ by a plane passing through the other two vertices of $\Delta$ and cutting through all other edges incident to $L(v)$ at a sufficiently small distance from $L(v)$ (see **Fig. 1.7**).

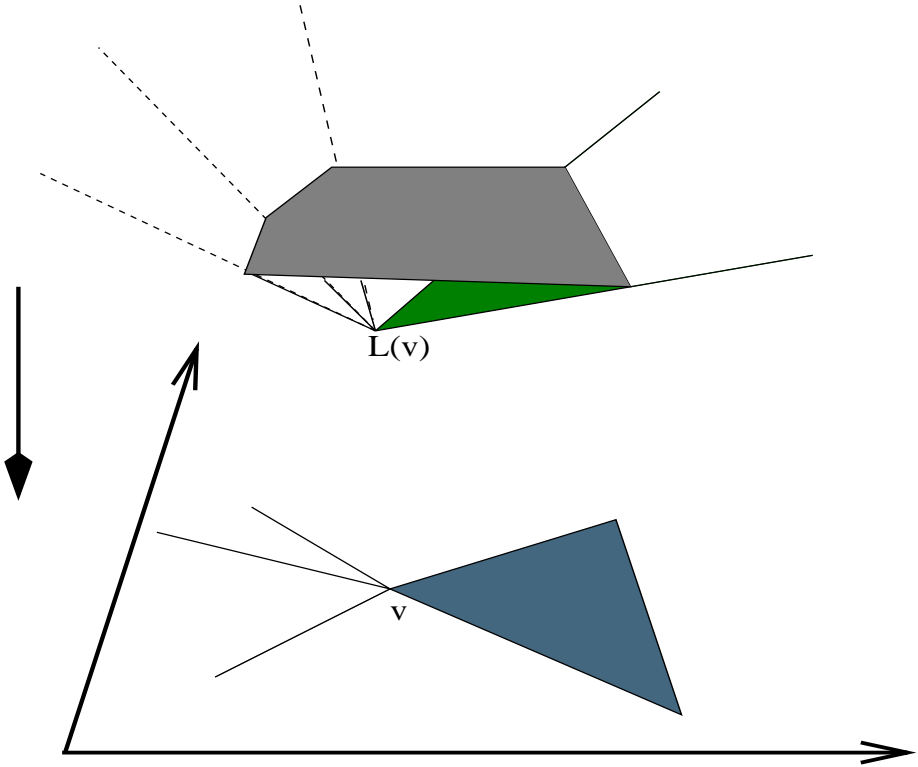

**Figure 1.7**

By Maxwell correspondence the projection of the resulting surface back on the plane gives a triangulation that supports tension in its 1-skeleton. Thus, if minimum subset exists then it consists of holes all whose vertices come from the vertices of $H$.

There is a finite number of way to construct a system of holes in $M$ that cover $H$, whose vertices come from vertices of $H$ and $M$, and whose complement supports tension. Pick among such systems of holes a system with minimal area. Denote it by $H_1$. We claim that this is, indeed, $H_{min}$. Assume that there is a proper subset $H_2$ of $H_1$ that satisfies the first two conditions.

The procedure described above can be repeatedly applied to $H_2$ until the vertex set of what is left of $H_2$ is a subset of the vertex set of $H$. But then the resulting collection of holes satisfies the first two conditions, and has smaller area than $H_1$. It contradict to the choice of $H_1$. Thus $H_1$ is minimal with respect to inclusion. Now, let's show that $H_1$ is actually contained in any system of holes satisfying the first two conditions. Suppose $H'$ is another system of holes which does not contain $H_1$, but satisfies the first two conditions. Consider a system of holes whose holes are intersections of the holes of $H'$ and $H_1$. Denote it by $H_{mesh}$. Since both $H'$ and $H_1$ covered $H$, $H_{mesh}$ covers $H$ too.

By definition $M\backslash H_1$ and $M\backslash H'$ both support tension, i.e. have partitions whose 1-skeletons support tension. Denote these graphs $G_1$ and $G'$. Consider a new graph, whose vertex set consists of all the vertices of the old graphs, $G_1$ and $G'$, and all intersections of edges of the old graphs ( this graph is what we see when $G_1$ and $G'$ are drawn on the same picture of $M$). The edge set of the new graph is defined accordingly. Denote it by $G_{mesh}$. Tensions on $G_1$ and $G'$ give rise to a tension on $G_{mesh}$. Thus, $G_{mesh}$ gives a partition of $M\backslash H_{mesh}$ whose 1-skeleton supports an all non-zero tension. Therefore we have a system of holes on $M$ that satisfies the first two conditions, and is contained in $H_1$. But we have already shown that this is impossible. Thus $H_1 = H_{min}$. □

If $\mathcal{H}$ is infinite we can still conclude that a maximal (with respect to inclusion) set supporting tension is the complement of a set of edge-disjoint polygons all whose vertices come from the vertex set of $H$, since the argument about cutting a corner can be also applied to an infinite surface.

As a consequence of the above theorem, the problem of finding the maximal area supporting tension can be rephrased as the following optimization problem on convex polyhedra. Find a convex surface in $\mathbb{R}^3$ such that each hole of $\mathcal{H}$ is covered by the vertical projection of a facet of the visible from the plane part of this surface, and the total area of the projections of the facets that cover holes is minimal possible. Another important implication of this theorem is that for finite systems of holes the complement of the set where tension *vanishes* must *support* tension, i.e. it has a triangulation which is a spider web. The case of infinite system of holes is more complicated. Even under additional restrictions on the system of holes, like that the vertices of the holes form an $(r, R)$-system, or that the sizes of the holes are uniformly bounded from both above and below it is not obvious that the union of all subsets of holes supporting tension can be represented as the complement of a discrete set of edge-disjoint polygons.

**Conjecture 1.12** *Let $\mathcal{H}$ be an infinite discrete system of polygons in $\mathbb{R}^2$. Then the union of all subsets of $\mathbb{R}^2 \backslash H$ supporting tension can be represented as the complement of a discrete set of edge-disjoint convex polygons.*

In addition we conjecture that the complement of the set where tension vanishes is not only polygonal, but also supports tension. Notice that our main result is about the *vanishing* of tension and the truth or false of this conjecture are irrelevant to our proof.

## 2 Tension Percolation Models

When not all of the complement of $H$ supports tension, it is natural to ask what is the maximum subset of $S\backslash H$ supporting tension.

Connelly and Rybnikov (1998, 1999) showed that the problem of finding the subset of maximum area supporting tension can be reduced to a linear programming problem, and as a consequence there is a polynomial (in the number of vertices of the holes) algorithm finding the maximal tensed subset. This polynomial algorithm can be effectively implemented and used for computer simulation of the membrane model for different distributions of holes and their parameters. However, this empirical evaluation for finite systems does not allow us to understand the limiting behavior of the system when the number of holes is infinite. In this section we consider a number of probability models for the distribution of holes and their sizes.

In the previous sections we looked at the behavior of tension on finite membrane with polygonal defects. A natural question is how the tension is redistributed on an infinite membrane (film) with polygonal holes. One can assume different distributions of the holes and their shapes.

The naïve model would be that all defects are centered at nodes distributed according to Poisson law, and the size of the hole is independent of the Poisson distribution. For this model we give a complete analysis and show that assuming that the holes are distributed according to Poisson law with parameter $\lambda$, and the shapes of the holes are IID and independent of the Poisson distribution, the tension ceases to exists on all of $\mathbb{R}^2$ for any value of Poisson parameter $\lambda$.

A more realistic approach would be as follows. Consider a Voronoi partition of the plane whose set of sites is a realization of a Poisson distribution with parameter $\lambda$. Declare each Voronoi domain normal with probability $p$ and defective with probability $1 - p$ (for percolation on Voronoi partitions see Vahidi-Asl and Wierman (1990)). For each defective domain define a convex hole centered at the Voronoi center of this domain. We suggest a family of model in which each convex polygonal hole is situated at a small neighborhood of its Voronoi domain. We associate with each hole three distributions: the distribution of the number of vertices of the hole, the distribution of the angular coordinate of a vertex of the hole, and that of the radial coordinate of the vertex. Let $\nu$ be any discrete distribution for the number of vertices of a defect. Let $\rho$ be a distribution for the length of the radius-vector of vertex. A reasonable choice for $\rho$ would be some random variable taking on all values between the radius of the inscribed circle and the radius of circumscribing circle of the domain. A natural choice for the distribution of the angular coordinate is a uniform distribution on the circle. The basic question is whether there is a critical value of parameter $p$ such that there is an infinite connected component of $\mathbb{R}^2$ supporting tension. A special case of such model would be one where the holes coincide with their Voronoi domains. We conjecture that in this case for any value of $\lambda$ and for any value of $p$ there is no infinite connected component supporting tension. The constructions involving Voronoi partitions are motivated by our desire to have a model in which the effect of *pinwheel configuration* is significant. In the "naïve" Poisson model the vanishing of tension is caused by the overlapping of holes, which is in turn is caused by the independence of the distribution of the hole sizes and the underlining Poisson process.

A lattice analog of continuous tension percolation models was introduced by Connelly and Rybnikov (1998, 1999). Let $\mathbb{T}$ be a triangular lattice on the plane. Remove each edge with probability $1 - p$. The question is whether there is a critical value of probability $p_c^t$ such when if $p$ is greater than $p_c^t$ there is an infinite subgraph of $\mathbb{T}$ supporting tension. This looks reasonable for modeling tension in a membrane, but the assumption that all atoms lie on three families of parallel lines does not seem to be realistic for describing properties of molecular systems like glasses. Thus, a generic version of this problem should have instead of $\mathbb{T}$ a generic triangular lattice. Note that the problem of tension percolation on a lattice is related to rigidity percolation on lattices. See Holroyd (1998) for rigorous estimates of the critical probability of rigidity percolation and Duxbury, Jacobs, Thorpe, Moukarzel (1995, 1996, 1997, 1998, 1999) for simulation results and their physical interpretation.

# 3  Percolation on iterative defects

In this section we show that if the centers of holes are distributed in $\mathbb{R}^2$ according to Poisson law and their shapes are IID, tension disappears on all of $\mathbb{R}^2$. Moreover, the following results on iterative convex hulls are true not only for dimension two, but also for $\mathbb{R}^d$. The proofs can be adopted straightforwardly. It is convenient for as to adopt a more general definition of hole, namely without the polygonality restriction.

**Definition 3.1** *A hole (f-hole) centered at $p \in \mathbb{R}^d$ is a region*

$$H(p, f) = \{p + f(\frac{\mathbf{x}}{\|\mathbf{x}\|})\mathbf{x} \mid \|\mathbf{x}\| \leq 1\}$$

*where $f$ is a continuous positive function defined on a unit $(d-1)$-sphere.*

Let $\mathcal{H}$ be some collection of holes.

**Definition 3.2** *We say that a collection of holes $\mathcal{H} = \{H_1(p_1, f_1), H_2(p_2, f_2), \ldots\}$ is uniformly bounded from below, if there exist $r > 0$ such that $f_i(\mathbf{x}) \geq r$ for all $f_i$ and $\|\mathbf{x}\| = 1$. Similarly, we say that $\mathcal{H}$ is bounded from above, if there exist $R > 0$ such that $f_i(\mathbf{x}) \leq r$ for all $f_i$ and $\|\mathbf{x}\| = 1$.*

Consider a $d$−dimensional Poisson point process with rate $\lambda$. Let $Y = Y(\omega)$ be the collection of nodes of some realization $\omega$ of the process. Each node $y \in Y(\omega)$ is the center of a hole $H(y, f_y)$, where function $f_y$ is positive and continuous.

Let $\mu$ be a probability measure on some subspace of $C(\mathbb{S}^{d-1}, \mathbb{R}_+)$. Suppose that for each $y$ the function $f_y$ is chosen from a distribution $\mu$ independently of the other functions and the configuration $\omega$. Therefore, the holes $H(y, f_y)$ are IID. Under this condition, it is known that if $\mathcal{H}$ is bounded both from above and below, there exists some $\lambda_c = \lambda_c(\mathcal{H}, \mu)$ such that for $\lambda < \lambda_c$ there exists no infinite cluster and for $\lambda > \lambda_c$ such cluster exists a.s (see Menshikov and Sidorenko (1987)).

The next statement is our main percolation result. In part, it can be viewed as a generalization of some results of Schonmann (1992) for bootstrap percolation on $\mathbb{Z}^d$.

**Theorem 3.3** *For any distribution $\mu$ and any $\lambda > 0$ there exists a nonnegative integer $N = N(\mu, \lambda)$ such that a.s. there is an infinite $N$-cluster whose convex hull coincides with $\mathbb{R}^d$.*

To prove this theorem, we reduce it in two steps (see Lemmas 3.4 and 3.5) to a simpler proposition (see Proposition 3.7).

**Lemma 3.4** *It is sufficient to prove Theorem 3.3 only for the case when the set of defects is bounded from below by some $r > 0$.*

*Proof.* According to definition, a hole corresponds to a continuous
and positive function $f(u)$ on a unit sphere $\mathbb{S}^{d-1}$. Since $\mathbb{S}^{d-1}$ is compact and $f(u)$ is continuous, $f(u)$ achieves its minimum there, and this minimum is strictly positive. Let $\rho(f) := \min_{u \in U} f(u) > 0$. The distribution $\mu$ on functions $f$ induces a distribution on a real-valued *positive* random variable $\rho$. Therefore, there must exist $r > 0$ such that $P(\rho \geq r) = p_0 > 0$.

Now we couple the original process with rate $\lambda$ and the functions $f_y$ distributed according to $\mu$, with a process which culls all the defects for which $\rho(f) < r$. It is straightforward that this process is also a Poisson process, but with rate $\tilde{\lambda} = \lambda p_0$. The distribution of defects $\tilde{\mu}$ for the thinned process is such that the radius of the inscribed circle for *any* defect is at least $r$. Consequently, the set of holes is now uniformly bounded below. If we show that the statement of Theorem 3.3 holds for $\tilde{\mu}$ and any $\tilde{\lambda} > 0$, it will imply Theorem 3.3 for a generic $\mu$, as the set of holes in "tilde" model is stochastically smaller than that of the original one. □

Hence, we can assume that each hole ($\equiv$ defect of $0^{th}$ generation) corresponding to a node $y \in Y$ contains a ball of radius $r$ with the center at $y$. Since taking a convex hull is a monotonous operation (namely, if $A \subseteq B$ then $conv(A) \subseteq conv(B)$), it suffices to prove Theorem 3.3 for the case when all holes are the balls of radius $r$. Moreover, without loss of generality we will can assume that $r = 1$ as we can always re-scale the space. Consequently, Theorem 3.3 follows from the next lemma.

**Lemma 3.5** *Let $\mathcal{H}$ be a collection of balls of unit radius whose centers are distributed in $\mathbb{R}^d$ according to a Poisson law with parameter $\lambda > 0$. Then there exists a nonnegative integer $N = N(\lambda)$ such that a.s. there exists an infinite $N$-cluster whose convex hull is $\mathbb{R}^d$.*

## 3.1 Proof of Lemma 3.5

We will present the proof for the case $d = 2$ (so that the defects are in $\mathbb{R}^2$). The arguments can be generalized for $d \geq 3$ rather easily.

Pick a hole(which is a unit circle) from $\mathcal{H}$, and denote it by $G(0)$. Examine all circles from $\mathcal{H}$ which have non-empty intersection with $G(0)$, take the convex hull of them and $G(0)$, and denote the resulting set by $G(1)$. Then consider all circles of $\mathcal{H}$ which intersect $G(1)$, and call the convex hull of $G(1)$ and their union $G(2)$. Iterating this procedure we construct the sets $G(2), G(3), \ldots$. If for some $k$ no element of $\mathcal{H}$ intersects with $G(k)$, we have $G(k+l) \equiv G(k)$ for all $l \geq 0$.

**Proposition 3.6** *$G(k)$ is contained in some defect of $k^{th}$ generation.*

This statement is obvious, because the operation of taking convex hull is commutative.

For any point $a \in \mathbb{R}^2$ let $(\rho_a, \varphi_a)$ be its polar and $(x_a, y_a)$ be its Cartesian coordinates. Without loss of generality assume that $G(0)$ is centered at the origin. Denote by $\mathcal{C}(k)$ the circle of radius $k$ centered at the origin.

**Proposition 3.7** *With a positive probability depending on $\lambda$ only,*

$$\mathcal{C}(k+1) \subseteq G(k) \text{ for all } k = 1, 2, 3, \ldots \tag{3}$$

*that is, $G(k)$'s eventually cover all $\mathbb{R}^2$ as $k \to \infty$.*

*Proof.* Here we actually prove a stronger result: we show that (3) holds even if $G(k+1)$ were constructed using $G(k)$ and *only* those circles of $\mathcal{H}$ whose centers lie inside the ring

$$R_{k+1} = \{a \in \mathbb{R}^2 : k + \frac{3}{2} \leq \rho_a \leq k + 2\}.$$

$G(\cdot)$ is an increasing process taking values in the subsets of $\mathbb{R}^2$. Observe that $G(k-1) \subseteq \mathcal{C}(k+1)$), and let $E_k$ be the event "$\mathcal{C}(k+1) \subseteq G(k)$". Observe that the "smallest" value $G(k)$ can take on $E_k$ is $\mathcal{C}(k+1)$. Therefore, if we can demonstrate that

$$\mathsf{P}(E_{k+1} \mid G(k) = \mathcal{C}(k+1), G(k-1), G(k-2), \ldots, G(1))$$
$$= \mathsf{P}\left(E_{k+1} \mid G(k) = \mathcal{C}(k+1)\right) \geq 1 - \gamma_k, \tag{4}$$

it will imply that $\mathsf{P}(E_{k+1} \,|\, E_k, E_{k-1}, \ldots, E_1) \geq 1 - \gamma_k$ and as a result

$$\mathsf{P}(\bigcap_{k=1}^{\infty} E_k) = \prod_{k=1}^{\infty}(1 - \gamma_k) > 0 \tag{5}$$

as soon as $\sum_k \gamma_k < \infty$.

It is straightforward that $\mathsf{P}(E_{k+1} \,|\, G(k) = \mathcal{C}(k+1)) > 0$ for all $k$, so it suffices to estimate $\gamma_k$ from (4) only for large values of $k$. Suppose that $G(k) = \mathcal{C}(k+1)$ and break the ring $R_k$ into $M = M(k) := 2\pi\sqrt{k+3}$ congruent pieces

$$R_{k+1}^{(i)} = \{a \in R_{k+1} : \varphi_a \in (\alpha i, \alpha(i+1)]\}, \quad i = 0, 1, \ldots, M-1$$

where $\alpha = \alpha(k) = \dfrac{2\pi}{M} = \dfrac{1}{\sqrt{k+3}}$

Let the event $E_{k+1}^{(i)}$ be "$R_{k+1}^{(i)}$ contains at least one center of a circle from $\mathcal{H}$". Then these events are independent for different values of $i$ and $k$, and consequently the probability of the event $\tilde{E}_{k+1} = \cap_{i=0}^{M-1} E_{k+1}^{(i)}$ "each $R_{k+1}^{(i)}$ contains such center" is

$$\mathsf{P}(\tilde{E}_{k+1}) = \left(1 - e^{-\lambda\alpha(k+7/4)/2}\right)^M \geq 1 - Me^{-\lambda\alpha(k+7/4)/2} = 1 - C_1 k^{1/2} e^{-C_2 k^{1/2}}$$

for some positive constants $C_1$ and $C_2$ and large $k$.

Our next step is to show that $\tilde{E}_{k+1}$ implies $E_{k+1}$ whenever $\mathcal{C}(k+1) \subseteq G(k)$. Once this is established, we have

$$\gamma_k \leq C_1 k^{1/2} e^{-C_2 k^{1/2}},$$

and, as a result, $\sum \gamma_k < \infty$, which yields (5) and proves the proposition.

Assume that $\tilde{E}_{k+1}$ occurs and consider an arbitrary point

$$a \in \mathcal{C}(k+2). \tag{6}$$

Since each of $R_{k+1}^{(i)}$ contains a center of a circle of $\mathcal{H}$ and every such circle obviously intersects $\mathcal{C}(k+1)$, there exist two points $b$ and $c$ in $G(k+1)$ such that

$$\rho_b \geq k + 5/2, \qquad \rho_c \geq k + 5/2,$$
$$\varphi_b \leq \varphi_a \leq \varphi_c \mod 2\pi,$$
$$\varphi_a - \varphi_b \leq \alpha \mod 2\pi, \qquad \varphi_c - \varphi_a \leq \alpha \mod 2\pi. \tag{7}$$

Since $G(k+1)$ is convex, it must contain all points inside the triangle formed by the origin $0$ and points $b$ and $c$. For any two points $e$ and $f$ of the plane, let $[e, f]$ denote the segment with these points being its endpoints. There exist a point $h$ lying on the segment $[b, c]$ such that $[0, h]$ is orthogonal to $[b, c]$. Notice that this point has the smallest distance from the origin among all points of the segment $[b, c]$. Let us get the lower bound on $\rho_h$. From (7) it follows that at least one of the angles between $[0, h]$ and $[0, b]$ or $[0, c]$ does not exceed $\alpha$. Suppose it is the first one. The length of $[h, b]$ in this case is smaller or equal to

$$\rho_b \sin\alpha \leq \rho_b \alpha = \dfrac{\rho_b}{\sqrt{k+3}}$$

whence

$$\rho_h \geq \sqrt{\rho_b^2 - \frac{\rho_b^2}{k+3}} \geq \sqrt{\frac{(k+5/2)^2(k+2)}{k+3}} \geq k+2$$

as $(k+5/2)^2 > (k+2)(k+3)$. Consequently, equations (7) and 6 yield that the point $a$ lies inside the triangle $\Delta(b0c)$ and hence belongs to $G(k+1)$. □

Let us make some observations which rely on obvious generalizations of Proposition 3.7. In our discussion we will rely on the proof presented above.

First, color the plane as an infinite chess-board with the cell size of $1/4$, so that a point with coordinates $(x,y)$ is green if $\lfloor 4x \rfloor + \lfloor 4y \rfloor$ is even, and blue if this sum is odd. The arguments in the proof will remain valid if instead of $\mathcal{H}$ we consider only those circles of it which have "green" centers (these centers form a Poisson point process with the same rate $\lambda$ on the "green part" of the plane). The similar statement about the "blue" process is valid as well.

Next, note that if while constructing $G(k+1)$ from $G(k)$ we do not use any circle with a center $a$ such that

$$\begin{aligned} -k^{-2/3} \leq \varphi_a \leq k^{-2/3} \text{ mod } 2\pi \text{ or} \\ \pi - k^{-2/3} \leq \varphi_a \leq \pi + k^{-2/3} \end{aligned} \quad (8)$$

(and use only the green part of the plane), the arguments above still hold for large $k$, as $k^{-2/3} = o(\alpha(k))$. Denote the resulting process $G'(k)$, and the probabilities corresponding to (4) as $1-\gamma'_k$. Since $\gamma'_k$ are asymptotically of the same form as $\gamma_k$, we have $\sum \gamma'_k < \infty$.

The third observation is that if the circle $\mathcal{C}(k_0)$ $(k_0 > 1)$ is completely covered by some elements of $\mathcal{H}$ with all their centers lying inside $\mathcal{C}(k_0)$, then the probability $p(k_0)$ that the process $G'(k)$ described above (using only green part of the plane and "avoiding" points close to the horizontal axis), can be continued indefinitely, in a way described by (3), has the property

$$\lim_{k_0 \to \infty} p(k_0) = 1.$$

This follows from the fact that $p(k_0) = \prod_{k=k_0}^{\infty}(1-\gamma'_k) \geq 1 - \sum_{k \geq k_0} \gamma'_k$ and the sum on the right hand side converges to zero as $k_0 \to \infty$.

Consequently, there exist $k_0 > 1$ such that $p(k_0) \geq 0.9$. Fix $k_0$ and let $\nu$ denote the probability that indeed $\mathcal{C}(k_0)$ is completely covered by circles of $\mathcal{H}$ as described in the previous paragraph. The value $\nu$ may be small, however it is positive and therefore there exist a positive integer $T$ such that

$$(1-\nu)^{T+1} < 0.1 \quad (9)$$

Consider $T+1$ circles of radius $k_0$ on the plane with their centers located at the points

$$(0,0), \ (k_0^3, 0), \ (2k_0^3, 0), \ \ldots, \ (Tk_0^3, 0). \quad (10)$$

According to (9), the probability that at least one of them is completely covered by circles of $\mathcal{H}$ is larger than 0.9. Arbitrarily pick one such circle (let it be, say, the $i^{th}$ one). The event "the process $G'(k)$ starting from $(ik_0^3, 0)$ can be continued indefinitely" is independent of the configuration inside the other circles, since when the process $G'(\cdot)$ reaches the $j^{th}$ circle $(j \neq i)$ $k \approx |i-j|k_0^3 =: mk_0^3$, $m \geq 1$, and therefore the $j^{th}$ circle lies entirely inside the angle

$$-\frac{1}{mk_0^2} < \varphi < \frac{1}{mk_0^2} \text{ or } \pi - \frac{1}{mk_0^2} < \varphi < \pi + \frac{1}{mk_0^2}$$

and at the same time

$$\frac{1}{mk_0^2} \leq \frac{1}{m^{2/3}k_0^2} \approx \frac{1}{k^{2/3}}.$$

Hence, with probability exceeding $(1 - (1-\nu)^{T+1}) \times p(k_0) \geq 0.81$ one of the $k_0$-radius circles with a center at one of the points of (10) is completely covered by elements of $\mathcal{H}$ *and* the process $G'(k)$ which starts from $i^{th}$ circle can grow indefinitely.

Now let us tile the plane with the boxes $L \times L$ where $L = 5Tk_0^3$, so that for any $(i,j) \in \mathbb{Z}^2$, the box $B(i,j)$ consists of the points

$$\{(iL + x', jL + y') \text{ where } (x', y') \in (0, L]^2\}$$

Suppose that $i+j$ is an even number. We say that $B(i,j)$ is *open* (see **Fig. 3.8**), if

(A) one of $T+1$ circles of radius $k_0$ with the centers in the set

$$(X_i, Y_j), \ (X_i + k_0^3, Y_j), \ (X_i + 2k_0^3, Y_j), \ \ldots, \quad (X_i + Tk_0^3, Y_j)$$
$$\equiv ((i+0.6)L, Y_j)$$
$$\text{where } X_i = (i+0.4)L, \ Y_j = (j+0.5)L$$

is completely covered by the elements of $\mathcal{H}$ *and*

(B) the process $G'(\cdot)$ starting with the center $a_{ij}$ of a circle satisfying (A) and using only green points for $k > k_0$ reaches the circumference of radius of $0.6L$ (i.e. $G'(k_1 - 1)$, where $k_1 = \lceil 0.6L \rceil$, contains the entire circle $a_{ij} + \mathcal{C}(k_1)$).

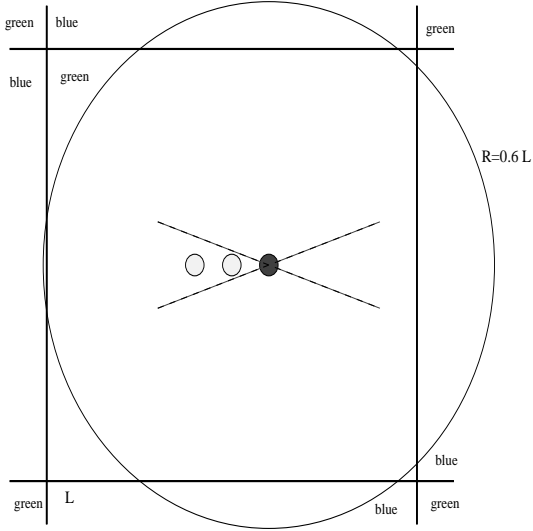

**Figure 3.8**

In the case when $i+j$ is odd, we call $B(i,j)$ open if the same is true except $G'(k)$ uses *blue* points.

From the translation invariance and the previous arguments it follows that

$$\mathsf{P}(B(i,j) \text{ is open}) > 0.8$$

Moreover, the boxes $B(i,j)$ are open independently of each other, since the neighboring (by an edge) boxes use points of different colors and $G'(k_1 - 1)$ does not reach any of the corners of the box (since $0.6L < \sqrt{(0.5L)^2 + (0.4L)^2}$). Also observe that when two boxes $B(i,j)$ and $B(i',j')$ sharing a common edge are open, the circles $\mathcal{C}(k_1) + a_{ij}$ and $\mathcal{C}(k_1) + a_{i'j'}$ have a non-empty intersection.

We say that two boxes are connected, if there exists a sequence of open boxes starting with the first one and ending with the second one such that each pair of two consecutive boxes in this sequence share an edge. The "box"-cluster containing the box $B(i,j)$ is the set of all boxes connected to it.

Hence, we can couple the boxes $\{B(i,j)\}$ with the site percolation process on $\mathbb{Z}^2$ where each site is open independently of the others with probability exceeding 0.8. However, as follows from Zuev (1987), the critical probability $p_c$ for the site percolation on $\mathbb{Z}^2$ is smaller than $0.68 < 0.8$ and therefore there exists a.s. an infinite open cluster of boxes $B(i,j)$. As a result, Proposition 3.6 yields that there exist a.s. an infinite $(k_1 - 1)$-cluster.

Let us prove the following

**Proposition 3.8** *Consider a site percolation model on $\mathbb{Z}^2$ in the supercritical regime (i.e. $p > p_c$). Then the convex hull of a (unique) infinite cluster is $\mathbb{R}^2$.*

*Proof of the proposition.* First, we remark that the fact that the infinite cluster is unique whenever it exists follows from Aizenman, Kesten and Newman (1987). Let $C(\infty)$ be such infinite cluster. On the plane with coordinates $(x,y)$ consider the angle

$$A = \{(x,y) \in \mathbb{Z}^2 : x \geq 0, \ 0 \leq y \leq x\}.$$

Notice that the plane can be partitioned into 8 cones congruent to $A$. Let $E_A$ be the event $|A \bigcap C(\infty)| = \infty$, i.e. there are infinitely many points of the infinite cluster lying inside of $A$. Since this is a tail event, $\mathsf{P}(E_A)$ is either 0 or 1. If it were 0, it would imply that the probability that there are infinitely many points of $C(\infty)$ in each of the other seven congruent to $A$ cones is also 0. This contradicts to the fact that $|C(\infty)| = \infty$. Hence, $P(E_A) = 1$ and the similar is true for the others cones. However, it is easy to see from geometrical observations that if each of this eight cones contains infinitely many points of $C(\infty)$ then $conv(C(\infty)) = \mathbb{R}^2$. □

Hence, the convex hull of the $(k_1 - 1)$-cluster is at least as large as that of an infinite open cluster for the site percolation on $\mathbb{Z}^2$ (rescaled $L$ times), and Lemma 3.5 has been proven. □

**Acknowledgments** K. R. thanks Robert Connelly for introducing him to the problem and many stimulating discussions. S. V. thanks Michael Aizenman and Roberto Schonmann for pointing out to him connections between the processes studied in the paper and bootstrap percolation models.